\documentclass[12pt]{amsart}
\usepackage{amsmath,amsthm,amssymb}
\usepackage{bm}
\usepackage{color}
\usepackage{graphicx}
\usepackage{xcolor}

\begin{document}
\newtheorem{theorem}{Theorem}[section]
\newtheorem{corollary}[theorem]{Corollary}%[section]
\newtheorem{lemma}[theorem]{Lemma}%[section]
\newtheorem{proposition}[theorem]{Proposition}%[section]
\theoremstyle{definition}
\newtheorem{definition}[theorem]{Definition}%[section]
\newtheorem{example}[theorem]{Example}%[section]
\newtheorem{remark}[theorem]{Remark}%[section]

\numberwithin{equation}{section}

\begin{center}
 Uniqueness and multiplicity of positive radial solutions 
 to the super-critical Brezis-Nirenberg problem in an annulus
\end{center}

\vspace{2ex}

\begin{center}
 Naoki Shioji \\[1ex]
 Department of Mathematics, Faculty of Engineering, \\
 Yokohama National University \\ 
 Tokiwadai, Hodogaya-ku, Yokohama 240--8501, Japan \\
 email: shioji193@gmail.com \\[2ex]
 Satoshi Tanaka 
 \\[1ex]
 Mathematical Institute, Tohoku University \\ 
 Aoba 6--3, Aramaki, Aoba-ku, Sendai 980--8578, Japan \\
 email: satoshi.tanaka.d4@tohoku.ac.jp \\
 and \\[2ex]
 Kohtaro Watanabe 
 \footnote{This work was supported by KAKENHI (18K03387). \\ 
 \hfill \today}
 \\[1ex]
 Department of Computer Science, National Defense Academy, \\
 1--10--20 Hashirimizu, Yokosuka 239--8686, Japan \\
 email: wata@nda.ac.jp
\end{center}

\bigskip

\noindent{\bf Abstract.}
The super-critical Brezis-Nirenberg problem in an annulus is considered.
The new uniqueness result of positive radial solutions is established
for the three-dimensional case.
It is also proved that the problem has at least three positive radial solutions 
when the inner radius of the annulus is sufficiently small and 
the outer radius of the annulus is in a certain range.
Moreover, for each positive integer $k$, the problem has at least $k$ positive 
radial solutions when the exponent of the equation is greater than
the critical Sobolev exponent and is less than the Joseph-Lundgren exponent.

\vspace{2ex}


\section{Introduction}

We consider the following Brezis-Nirenberg problem
\begin{equation}
 \left\{
  \begin{array}{cl}
   \Delta v + \lambda v + v^p = 0, & \ x\in A(a,b), \\[1ex]
    v=0, &  x\in \partial A(a,b),
  \end{array}
 \right. 
 \vspace{-1ex}
 \label{PDE}
\end{equation}
where $\lambda>0$, $p>1$, $A(a,b)=\{ x \in \mathbb{R}^N : a<|x|<b \}$,
$N\in\mathbb{N}$, $N\ge 2$, and $0<a<b<\infty$.
Every positive radial solution of \eqref{PDE} satisfies
\begin{equation}\label{A}
 \left\{
  \begin{array}{l}
   u'' + \dfrac{N-1}{r} u' + \lambda u + u^p = 0, \quad r \in (a,b), 
    \\[2ex]
   u(a)=u(b)=0.
  \end{array}
 \right.
\end{equation}
The following existence result of positive solutions to \eqref{A} is 
well-known and its proof will be given in Appendix for the reader's convenience.

\begin{proposition}\label{existpossol}
 Assume $\lambda>0$, $p>1$, $N\in\mathbb{N}$, $N\ge 2$ and $0<a<b<\infty$.
 Then problem \eqref{A} has at least one positive solution if and only if 
 $\lambda<\lambda_1(A(a,b))$.
\end{proposition}

Here and hereafter $\lambda_1(A(a,b))$ denotes the first eigenvalue of 
\begin{equation*}
 \left\{
  \begin{array}{cl}
   -\Delta \phi = \lambda \phi, & \ x\in A(a,b), \\[1ex]
    \phi=0, &  x\in \partial A(a,b).
  \end{array}
 \right. 
\end{equation*}

Uniqueness of solutions to \eqref{A} has been studied by
Ni and Nussbaum \cite{NN}, Yadava \cite{Yad}, Korman \cite{Kor}, 
Yao, Li and Chen \cite{YLC}, 
Shioji and Watanabe \cite{SW}, and 
Shioji, Tanaka and Watanabe \cite{STW}.

The following result is obtained by Yao, Li and Chen \cite[Theorem 1.1]{YLC}.

\bigskip

\noindent{\bf Theorem A (Yao, Li and Chen \cite[Theorem 1.1]{YLC}).}
\it
Let $N\in\mathbb{N}$, $N\ge2$, $0<a<b<\infty$ and 
$0<\lambda<\lambda_1(A(a,b))$.
Then problem \eqref{A} has a unique positive solution, if one of the following
conditions holds\textup{:}
\begin{enumerate}
 \item $N=2$, $p>1$\textup{;}
 \item $N\ge3$, $1<p\le(N+2)/(N-2)$\textup{;}
 \item $N\ge3$, $p>(N+2)/(N-2)$ and $r_0 \not\in(a,b)$, where
       \begin{align*}
	r_0=\sqrt{\frac{2((N-2)p+N-4)((N-2)p-(N+2))}{\lambda(p-1)(p+3)^2}}.
       \end{align*}
\end{enumerate}
\rm

\bigskip

The following result is established by Shioji and Watanabe \cite[Theorem 1]{SW}.

\bigskip

\noindent{\bf Theorem B (Shioji and Watanabe \cite[Theorem 1]{SW}).}
\it
Let $N=3$, $p>5$, $0<a<b<\infty$, $0<\lambda<\lambda_1(A(a,b))$ and 
$a+b\ge\pi/\sqrt{\lambda}$.
Then problem \eqref{A} has a unique positive solution.
\rm

\bigskip

Here, we would like to note that (original) Theorem 1 of \cite{SW} requires 
the additional assumption $a+b\ge\pi/\sqrt{\lambda}$.

We note that $p=5=(N+2)/(N-2)$ for $N=3$.
Since $\lambda_1(A(a,b))=\pi^2/(b-a)^2$ when $N=3$, 
the inequality $\lambda<\lambda_1(A(a,b))$ means
$b-a<\pi/\sqrt{\lambda}$.

The first main result of this paper is as follows.

\begin{theorem}\label{uniquenessresult}
 Let $N=3$, $p>5$ and $0<a<b<\infty$, and let $r_0$ be the constant as in 
 Theorem A.
 Assume either cases \textup{(i)} or \textup{(ii)} or \textup{(iii):}
 \begin{enumerate}
  \item $r_{0}\leq a \text{ or } b\leq r_{0}$\textup{;}
  \item $a<r_{0} \text{ and } b>r_{0}$ and
	\begin{align*}
	 a+\dfrac{\pi}{\sqrt{\lambda}}-\dfrac{1}{\sqrt{\lambda}}\arctan\left(\dfrac{2(p+3)\sqrt{\lambda}}{p-5}a\right)\le b<a+\dfrac{\pi}{\sqrt{\lambda}};
	\end{align*}
 \item $a<r_{0} \text{ and } b>r_{0}$ and
	\begin{align*}
	 b \le -a +\frac{1}{\sqrt{\lambda}}\left(\dfrac{2\sqrt{2(p-5)}}{p+3}+\arccos\left(-\sqrt{\dfrac{p-5}{p+3}}\right)\right).
	\end{align*}
 \end{enumerate}
 Then, problem \eqref{A} has a unique positive solution.
\end{theorem}
\begin{figure}[h]
  \centering
  \includegraphics[scale=0.72]{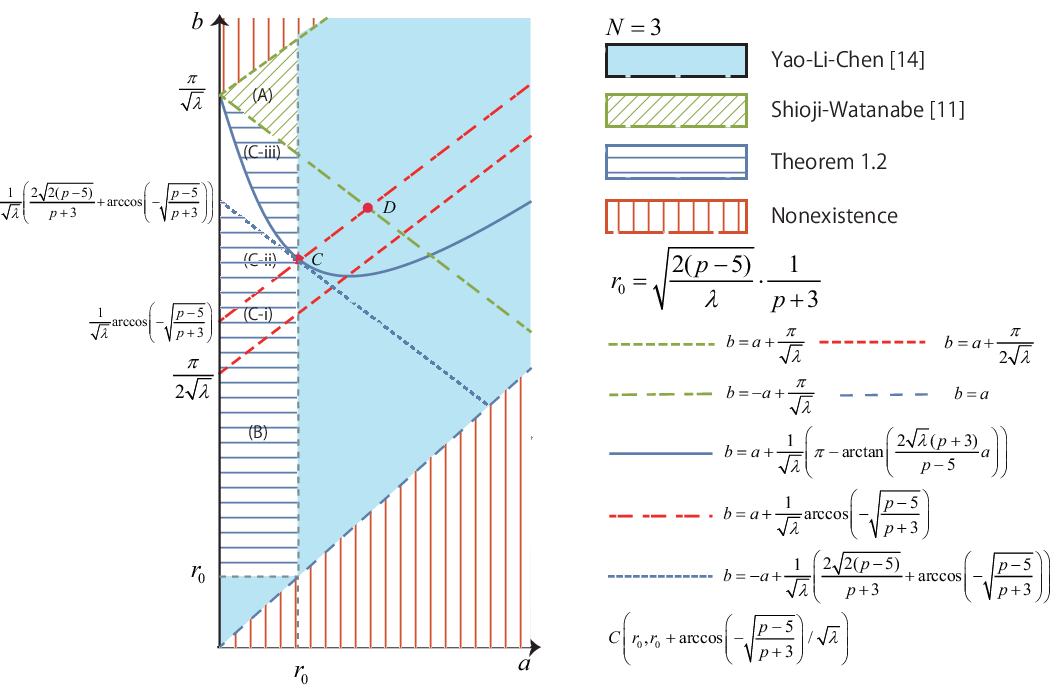}
  \caption{The regions where the uniqueness of solutions corresponding to Theorems A, B, and 1.2 is demonstrated when $N=3$ and $p>5$.}
  \label{fig:fig1}
\end{figure}
The case (i) of Theorem \ref{uniquenessresult} should be noted as the $N=3$ case of Theorem A.

Put $A(0,2\sqrt{2(p-5)}/(\sqrt{\lambda}(p+3))+\arccos(-\sqrt{(p-5)/(p+3)})/\sqrt{\lambda})$ and $B(0,\pi/\sqrt{\lambda})$.
Let $N=3$ and $p>5$. Moreover, let the domain enclosed by line segments $\overline{AB}$ and $\overline{CA}$ and arc $\stackrel{\frown}{BC}$ be denoted as $S_{p,\lambda}$. We note $S_{p,\lambda}$ is a region where the uniqueness of positive radial solutions of the problem \eqref{A} is not guaranteed. It should be noted that when $p$ is large, the measure of $S_{p,\lambda}$ becomes very small, since $r_{0}\rightarrow 0$ as $p\rightarrow\infty$ and Remark \ref{limpinfty}.
In spite of this, we can show that problem \eqref{A} has at least three positive solutions for some $(a,b)\in S_{p,\lambda}$.
The phenomenon of the multiple existence of solutions on such an extremely small possible parameter set can 
also be observed in \cite{STW2}.
\begin{remark}
We note that
\begin{equation*}
 b = a+\dfrac{\pi}{\sqrt{\lambda}}-\dfrac{1}{\sqrt{\lambda}}\arctan\left(\dfrac{2(p+3)\sqrt{\lambda}}{p-5}a\right)\quad\quad (=\varphi_{1}(a)),
\end{equation*}
\begin{equation*}
 b =  a + \dfrac{1}{\pi}\arccos\left(-\sqrt{\dfrac{p-5}{p+3}}\right)\quad\quad (=\varphi_{2}(a)),
\end{equation*}
and
\begin{equation*}
 b = -a +\frac{1}{\sqrt{\lambda}}\left(\dfrac{2\sqrt{2(p-5)}}{p+3}+\arccos\left(-\sqrt{\dfrac{p-5}{p+3}}\right)\right)\quad (=\varphi_{3}(a))
\end{equation*}
intersect at $a=r_0$ (hence, these three curves intersect at point C of Fig. \ref{fig:fig1}). Moreover, it holds
\begin{equation*}
 \dfrac{\partial^2 \varphi_{1}}{\partial a^2}(a)=\frac{16 a \lambda  (p-5) (p+3)^3}{\left(4 a^2 \lambda  (p+3)^2+(p-5)^2\right)^2}>0
\end{equation*}
and 
\begin{equation*}
\dfrac{\partial\varphi_{1}}{\partial a}(a)\Big|_{a=r_0}=-1,
\end{equation*}
Therefore, the positional relationship of each graph will be as shown in Fig. \ref{fig:fig1}.
\end{remark}

\begin{theorem}\label{exist3sols}
 Let $N\in\mathbb{N}$, $N\ge3$, $p>(N+2)/(N-2)$, $b>0$ and $\lambda>0$.
 If the problem on a ball
 \begin{equation}
 \left\{
  \begin{array}{ll}
   u'' + \dfrac{N-1}{r} u' + \lambda u + u^p = 0, &  r\in(0,b), \\[2ex]
   u'(0)=u(b)=0, & 
  \end{array}
 \right. 
 \label{B}
 \end{equation}
 has at least one positive solution, then there exists $a_0\in(0,b)$ 
 such that problem \eqref{A} has at least three positive solutions for 
 $0<a<a_0$.
\end{theorem}

There are a lot of study on the ball.
First, Brezis and Nirenberg \cite{BN} considered the critical case 
$p=(N+2)/(N-2)$.
The super-critical case $p>(N+2)/(N-2)$ has been studied by
Merle and Peletier \cite{MP},
Dolbeault and Flores \cite{DF},
Guo and Wei \cite{GW}, 
Miyamoto \cite{M2018},
and
Miyamoto and Naito \cite{MN2020, MN2023, MN2024}.

The following result has been obtained essential by
Dolbeault-Flores \cite[Theorem 1]{DF}, 
Guo-Wei \cite[Theorems 1.1]{GW} and 
Miyamoto-Naito \cite[Corollary 1.1 and Remark 1.3 (i)]{MN2020}.
We will give a proof of Theorem C in Appendix for the reader's convenience.

\medskip

\noindent{\bf Theorem C.}
\it
Let $N\in\mathbb{N}$ and $N\ge3$.
Then there exist constants $\underline{r}>0$ and $r_*>0$ satisfying
the following {\rm (i)} and {\rm (ii):}  
\begin{enumerate}
 \item when $(N+2)/(N-2)<p<p_{JL}$, then problem \eqref{B} has at least one 
       positive solution for 
       $\underline{r}\le b<\sqrt{\lambda_1(B(1))/\lambda}$,
       and moreover, for each $k\in\mathbb{N}$, 
       problem \eqref{B} has at least $k$ positive solutions if
       $|b-r_*|$ is sufficiently small{\rm ;}
 \item when $p\ge p_{JL}$, then problem \eqref{B} has a positive 
       solution for $r_*<b<\sqrt{\lambda_1(B(1))/\lambda}$ 
\end{enumerate}
where $p_{JL}$ denotes the Joseph-Lundgren exponent  
\begin{equation*}
 p_{JL}:=\left\{ 
  \begin{array}{cl}
   1+\dfrac{4}{N-4-2\sqrt{N-1}}, & N \ge 11, \\[2ex]
   \infty, & 3\le N \le 10,
  \end{array}
 \right.
\end{equation*}
$B(b):=\{ x \in \mathbb{R}^N : |x|<b \}$
and $\lambda_1(B(b))$ is the first eigenvalue of 
\begin{equation*}
 \left\{
  \begin{array}{cl}
   -\Delta \phi = \lambda \phi, & \ x\in B(b), \\[1ex]
    \phi=0, &  x\in \partial B(b).
  \end{array}
 \right. 
\end{equation*}
\rm

\begin{remark}
 (i) It is well-known that the equation
 \begin{equation*}
  u'' + \dfrac{N-1}{r} u' + \lambda u + u^p = 0, \quad r>0, 
 \end{equation*}
 has a unique singular solution $u_*$, and $u_*$ satisfies 
 $\lim_{r\to0^+}u(r)=\infty$ and has a zero in $(0,\infty)$.
 The constant $r_*$ is the smallest zero of $u_*$ in $(0,\infty)$.
 See Section 3 and the proof of Theorem C in Appendix.

 (ii) It is known that if $\lambda\ge\lambda_1(B(b))$, then problem \eqref{B} 
 has no positive solution. 
 See Proposition A.2 in Appendix.
 Since $\lambda_1(B(b))=b^{-2}\lambda_1(B(1))$, we find that 
 $\lambda\ge\lambda_1(B(b))$ is equivalent to 
 $b\ge\sqrt{\lambda_1(B(1))/\lambda}$.

\end{remark}

Now we consider problem \eqref{A}.
Combining Theorems \ref{exist3sols} and C, we have the following corollary 
immediately.

\begin{corollary}
 Let $N\in\mathbb{N}$, $N\ge3$, $p>(N+2)/(N-2)$, $0<a<b<\infty$.
 Then there exists $a_{0}\in(0,b)$ such that problem \eqref{A} has at least 
 three positive solutions when $0<a<a_{0}$ for some $a_{0}\in (0,b)$, 
 if either following \textup{(i)} or \textup{(ii)} holds\textup{:}
 \begin{enumerate}
  \item $(N+2)/(N-2)<p<p_{JL}$ and 
	$\underline{r}\le b<\sqrt{\lambda_1(B(1))/\lambda}$\textup{;}
  \item $p\ge p_{JL}$ and $r_*<b<\sqrt{\lambda_1(B(1))/\lambda}$,
 \end{enumerate}
 where $\underline{r}$ and $r_*$ are positive constants as in Theorem C.
\end{corollary}

We also obtain the following result.

\begin{theorem}\label{existksols}
 Let $N\in\mathbb{N}$, $N\ge3$, $(N+2)/(N-2)<p<p_{JL}$, $\lambda>0$, and 
 let $r_*$ be a constants as in Theorem C.
 Then, for each $k\in\mathbb{N}$, there exists $\rho_k\in(0,r_*)$ such that
 if $|b-r_*|<\rho_k$, then problem \eqref{A} has at least $k$ positive 
 solutions when $0<a<a_0$ for some $a_0\in(0,b)$.
\end{theorem}

From Theorems \ref{uniquenessresult} and \ref{exist3sols}, 
we obtain the following corollary, which will be proved in Section 4.

\begin{corollary}\label{<lam<}
 Let $N=3$ and $p>5$. If \eqref{B} has a positive solution, then it holds
 \begin{equation*}
  \dfrac{1}{\sqrt{\lambda}}\left(\frac{2\sqrt{2(p-5)}}{p+3}+\arccos\left(-\sqrt{\dfrac{p-5}{p+3}}\right)\right)
  <b<\dfrac{\pi}{\sqrt{\lambda}}
 \end{equation*}
or equivalently
 \begin{equation}\label{Merle}
 \frac{1}{\pi^2}\cdot\left(\frac{2\sqrt{2(p-5)}}{p+3}+\arccos\left(-\sqrt{\dfrac{p-5}{p+3}}\right)\right)^{2}
  \cdot\lambda_{1}(B(b))
  <\lambda<\lambda_{1}(B(b)).
 \end{equation}
\end{corollary}

\begin{remark}\label{limpinfty}
 Letting $p\to\infty$, we have
 \begin{equation*}
  \left(\frac{2\sqrt{2(p-5)}}{p+3}+\arccos\left(-\sqrt{\dfrac{p-5}{p+3}}\right)\right)^{2} \to  \pi^2.
 \end{equation*}
Assume $N\geq 3$ and $p>(N+2)/(N-2)$.
We note Theorem A in Merle, Peletier and Serrin \cite{MPS} claims that if \eqref{B} has a positive solution, then $\lambda (=\lambda(p)) \rightarrow\lambda_{1}(B(b))$ as $p\rightarrow \infty$.
Therefore, Eq. \eqref{Merle} can be regarded as a refinement of Theorem A in \cite{MPS} for the case $N=3$.
\end{remark}

In Sections 2 and 3, we prove Theorem \ref{exist3sols} and \ref{existksols},
respectively.
We prove Theorem \ref{uniquenessresult} and Corollary \ref{<lam<} in Section 4. 
In Appendix, we give proofs of Proposition \ref{existpossol} and Theorem C.

\section{Proof of Theorem \ref{exist3sols}}

We prove Theorem \ref{exist3sols}.
The proof is based on the method of the proof of Theorem 1.10
in Ni and Nussbaum \cite{NN}, but our approach is simpler and more direct.

We assume that
$N\in\mathbb{N}$, $N\ge3$, $p>(N+2)/(N-2)$, $b>0$, $\lambda>0$, and 
problem \eqref{B} has a positive solution $u_0$.
Let $\phi$ be a unique solution of the initial value problem
\begin{equation}\label{phi}
 \left\{
  \begin{array}{l}
   \phi'' + \dfrac{N-1}{r} \phi' + \lambda \phi = 0, \quad r>0, \\[2ex]
   \phi(0)=1, \ \phi'(0)=0.
  \end{array}
 \right.
\end{equation}

\begin{lemma}\label{phi>0}
 $\phi(r)>0$ for $r\in[0,b]$.
\end{lemma}

\begin{proof}
 We note that $u_0$ and $\phi$ are solutions of 
 \begin{equation*}
  (r^{N-1} u_0')' + \lambda r^{N-1} (1 + \lambda^{-1}|u_0(r)|^{p-1}) u_0 = 0
 \end{equation*}
 and 
 \begin{equation*}
  (r^{N-1} \phi')' + \lambda r^{N-1} \phi = 0,
 \end{equation*}
 respectively.
 Assume to the contrary that $\phi$ has a zero in $(0,b]$, that is,
 $\phi(r)>0$ on $[0,r_0]$ and $\phi(r_0)=0$ for some $r_0\in(0,b]$.
 Then the Sturm comparison theorem leads us a contradiction.
 More precisely, 
 \begin{align*}
  r_0^{N-1} \phi'(r_0) u_0(r_0) 
   & = \int_0^{r_0} 
   \left[ r^{N-1} (\phi'(r)u_0(r)-\phi(r)u_0'(r) ) \right]' dr \\
   & = \int_0^{r_0} r^{N-1} |u_0(r)|^{p-1} u_0(r) \phi(r) dr > 0,
 \end{align*}
 which contradicts the fact that $\phi'(r_0)<0$ and $u_0(r_0)\ge0$.
\end{proof}

The equation 
\begin{equation*}
 u'' + \dfrac{N-1}{r} u' + \lambda u + |u|^{p-1} u = 0, \quad 0<r\le b 
\end{equation*}
is transformed into
\begin{equation}
 \frac{d^2U}{dt^2} + (N-2)^{-2} t^{-2(N-1)/(N-2)} (\lambda U + |U|^{p-1}U) = 0,
 \quad b^{-(N-2)} \le t < \infty,
 \label{EqU}
\end{equation}
 by 
\begin{equation*}
 U(t) = u(t^{-1/(N-2)}).
\end{equation*}
We define the functions $U_0$ and $\Phi$ by
\begin{equation*}
 U_0(t) = u_0(t^{-1/(N-2)}) \quad \textup{and} \quad 
 \Phi(t) = \phi(t^{-1/(N-2)}),
\end{equation*}
respectively. 
Then $U_0$ is a solution of \eqref{EqU} and $\Phi$ is a solution of
\begin{equation*}
 \frac{d^2\Phi}{dt^2} + \lambda (N-2)^{-2} t^{-2(N-1)/(N-2)} \Phi = 0,
 \quad b^{-(N-2)} \le t < \infty,
\end{equation*}
and we have
\begin{gather*}
  U_0(b^{-(N-2)})=0, \quad U_0(t)>0 \ \textup{for} \ t>b^{-(N-2)}, \\
   \lim_{t\to\infty} U_0(t) = u_0(0) \in (0,\infty), \\
  \Phi(t)>0 \ \textup{for} \ t \ge b^{-(N-2)}, \quad 
   \lim_{t\to\infty} \Phi(t) = 1.   
\end{gather*}
Moreover, we use the transformation
\begin{equation}\label{Utow}
 w(s) = \frac{U(t)}{\Phi(t)}, \quad 
 s = \int_{b^{-(N-2)}}^t \frac{1}{[\Phi(\tau)]^2} d\tau.
\end{equation}
Then \eqref{EqU} is transformed into
\begin{equation}
 \frac{d^2 w}{ds^2} + g(s)|w|^{p-1} w = 0,
 \quad s>0.
 \label{Eqw}
\end{equation}
where
\begin{equation}
 g(s) = (N-2)^{-2} t^{-2(N-1)/(N-2)}[\Phi(t)]^{p+3}, \quad
 s=\int_{b^{-(N-2)}}^t \frac{1}{[\Phi(\tau)]^2} d\tau. 
  \label{g}
\end{equation}
We note that $g(s)$ is positive and bounded on $[0,\infty)$. 
Now we define $w_0$ by
\begin{equation*}
 w_0(s) = \frac{U_0(t)}{\Phi(t)}, \quad 
 s=\int_{b^{-(N-2)}}^t \frac{1}{[\Phi(\tau)]^2} d\tau.
\end{equation*}
Then $w_0$ is a solution of \eqref{Eqw} and satisfies
$w_0(0)=0$, $w_0(s)>0$ for $s>0$ and 
$\lim_{s\to\infty} w_0(s)=u_0(0)\in(0,\infty)$.

We consider the following problem
\begin{equation}
 \left\{
  \begin{array}{l}
   w'' + h(s)|w|^{p-1} w = 0, \quad 0<s<c, \\[1ex]
   w(0)=w(c)=0, \\[1ex]
   w(s)>0, \quad 0<s<c,
  \end{array}
 \right.  
 \label{bvpw}
\end{equation}
where $h\in C[0,\infty)$, $h(s)>0$ for $s>0$, $p>1$ and 
$c>0$ is a fixed arbitrary constant. 
The following result will be shown later.

\begin{theorem}\label{exist3solsODE}
 Let $p>1$, $h\in C[0,\infty)$ and $h(s)>0$ for $s>0$.
 Assume that $h$ is bounded on $[0,\infty)$ and satisfy
 \begin{equation}
  \lim_{t\to\infty}t^{-(p+1)/2} \int_0^t s^{p+1} h(s) ds = \infty.
  \label{s^sigmah(s)>0}
 \end{equation}
 and that the equation
 \begin{equation}
  w'' + h(s)|w|^{p-1} w = 0
   \label{w''+hw^p=0}
 \end{equation}
 has a solution $w_0$ for which $w_0(0)=0$, $w_0(s)>0$ for $s>0$.
 Then problem \eqref{bvpw} has at least three solutions 
 for all sufficiently large $c>0$.
\end{theorem}

Now we give a proof of Theorem \ref{exist3sols}

\begin{proof}[Proof of Theorem \ref{exist3sols}]
 By L'Hopital's rule, we see that
 \begin{equation*}
  \lim_{s\to\infty} \frac{s}{t}
  = \lim_{t\to\infty} \frac{\int_{b^{-1/(N-1)}}^t \frac{1}{[\Phi(\tau)]^2} d\tau}{t}
  = \lim_{t\to\infty} \frac{1}{[\Phi(t)]^2}
  = 1.
 \end{equation*}
 Therefore,
 \begin{equation*}
  \lim_{s\to\infty} s^{2(N-1)/(N-2)} g(s) = (N-2)^{-2}.
 \end{equation*}
 We can check that \eqref{s^sigmah(s)>0} with $h=g$ is satisfied.
 Theorem \ref{exist3solsODE} implies that problem \eqref{bvpw} with $h=g$ 
 has at least three solutions for all sufficiently large $c>0$,
 which proves that Theorem \ref{exist3sols}.
\end{proof}

Let us prove Theorem \ref{exist3solsODE}.
We define the Rayleigh quotient
\begin{equation*}
 R(w):= \frac{\int_0^c |w'(s)|^2 ds}
  {\left(\int_0^c h(s) |w(s)|^{p+1} ds \right)^{\frac{2}{p+1}}}
\end{equation*}
and the energy
\begin{equation*}
 E(c) := \inf_{w \in H^1_0(0,c),\,w\not\equiv 0} R(w).
\end{equation*}
It is well-known that problem \eqref{bvpw} has a solution $w_c$ such that
\begin{equation*}
 R(w_c)=E(c).
\end{equation*}
We call such a solution $w_c$ a {\it least energy solution} of \eqref{bvpw}.

\begin{lemma}\label{E(c)-->0}
 Let $p>1$, $h\in C[0,\infty)$ and $h(s)>0$ for $s>0$.
 Assume that $h$ satisfy \eqref{s^sigmah(s)>0}.
 Then $E(c)\to0$ as $c\to\infty$.
\end{lemma}

\begin{proof}
 Let $c>4$.
 We define the function $W$ by
 \begin{equation*}
  W(s) = \left\{
  \begin{array}{ll}
   c^{-1/2} s, & 0 \le s \le c/2, \\[1ex]
   c^{1/2} - c^{-1/2} s, & c/2 \le s \le c.
  \end{array}
  \right.
 \end{equation*}
 Then $W\in C[0,c]$, $W(0)=W(c)=0$ and
 \begin{equation*}
  \int_0^c |W'(s)|^2 ds = 2 \int_0^{c/2} |W'(s)|^2 ds = 1.
 \end{equation*}
 We observe that 
 \begin{align*}
  \int_0^c h(s) |W(s)|^{p+1} ds 
   & \ge \int_{0}^{c/2} h(s) |W(s)|^{p+1} ds \\
   & = \int_{0}^{c/2} h(s) (c^{-1/2}s)^{p+1} ds \\
   & = 2^{(p+1)/2}(c/2)^{-(p+1)/2} \int_0^{c/2} s^{p+1} h(s) ds.
 \end{align*}
 Therefore,
 \begin{equation*}
  0 \le E(c) \le R(W)
    \le 2^{-1} 
    \left( (c/2)^{-(p+1)/2} \int_0^{c/2} s^{p+1} h(s) ds \right)^{-2/(p+1)}, 
 \end{equation*}
 which implies that $E(c)\to0$ as $c\to\infty$.
\end{proof}

\begin{lemma}\label{w'(0)-->0}
 Let $p>1$, $h\in C[0,\infty)$ and $h(s)>0$ for $s>0$.
 Assume that $h$ is bounded on $[0,\infty)$ and satisfy 
 \eqref{s^sigmah(s)>0}.
 Let $\{w_c\}_{c>0}$ be a family of least energy solutions to \eqref{bvpw}.
 Then $w_c'(0)\to0$ as $c\to\infty$. 
\end{lemma}

\begin{proof}
 Let $m\in(0,c)$ be a number satisfying $w_c'(m)=0$.
 Since 
 \begin{equation*}
  w_c''(s)=-h(s)(w_c(s))^p \le 0, \quad 0\le s \le c,
 \end{equation*}
 we see that $w_c'(0)\ge w_c'(s) \ge 0$ for $0\le s \le m$.
 Let $M>0$ be a constant $0\le h(s) \le M$ for $s\ge0$.
 By H\"{o}lder's inequality, we get
 \begin{align*}
  \frac{1}{2} (w_c'(0))^2 
   & = -\frac{1}{2} \int_0^m [(w_c'(s))^2]' ds \\ 
   & = -\int_0^m w_c'(s) w_c''(s) ds \\
   & = \int_0^m w_c'(s) h(s) (w_c(s))^p ds \\
   & \le \left( \int_0^m |w_c'(s)|^{p+1} ds \right)^{\frac{1}{p+1}}
         \left( \int_0^m |h(s)(w_c(s))^p|^{\frac{p+1}{p}} ds
              \right)^{\frac{p}{p+1}} \\
   & = \left( \int_0^m |w_c'(s)|^{p-1} |w_c'(s)|^2 ds 
       \right)^{\frac{1}{p+1}}
       \left( \int_0^m (h(s))^\frac{1}{p} h(s)|w_c(s)|^{p+1} ds
              \right)^{\frac{p}{p+1}} \\
   & \le (w_c'(0))^{\frac{p-1}{p+1}}
       \left( \int_0^m |w_c'(s)|^2 ds \right)^{\frac{1}{p+1}}
       M^\frac{1}{p+1} \left( \int_0^m h(s)|w_c(s)|^{p+1} ds
              \right)^{\frac{p}{p+1}}.
 \end{align*}
 Therefore,
 \begin{equation*}
  (w_c'(0))^\frac{p+3}{p+1} \le 2  M^\frac{1}{p+1}
  \left( \int_0^c |w_c'(s)|^2 ds \right)^{\frac{1}{p+1}}
  \left( \int_0^c h(s)|w_c(s)|^{p+1} ds
              \right)^{\frac{p}{p+1}}.
 \end{equation*}
 Multiplying the differential equation of \eqref{bvpw} by $w_c$ and 
 integrate it on $[0,c]$, we have
 \begin{equation*}
  \int_0^c |w_c'(s)|^2 ds = \int_0^c h(s) |w_c(s)|^{p+1} ds,
 \end{equation*}
 Hence, we find that
 \begin{equation*}
  (w_c'(0))^\frac{p+3}{p+1} \le 2 M^\frac{1}{p+1} \int_0^c |w_c'(s)|^2 ds
 \end{equation*}
 and
 \begin{equation*}
  E(c)= R(w_c) 
  = \frac{\int_0^c |w_c'(s)|^2 ds}
         {\left( \int_0^c |w_c'(s)|^2 ds \right)^{\frac{2}{p+1}}}
  = \left( \int_0^c |w_c'(s)|^2 ds \right)^{\frac{p-1}{p+1}},
 \end{equation*}
 which imply
 \begin{equation*}
  (w_c'(0))^\frac{p+3}{p+1} \le 2 M^\frac{1}{p+1} (E(c))^\frac{p+1}{p-1}.
 \end{equation*}
 From Lemma \ref{E(c)-->0}, it follows that $w_c'(0)\to0$ as $c\to\infty$.
\end{proof}

Let $w(s;\beta)$ denote a unique solution of initial value problem
\begin{equation*}
 \left\{
  \begin{array}{l}
   w'' + h(s)|w|^{p-1} w = 0, \quad 0<s<c, \\[1ex]
   w(0)=0, \quad w'(0)=\beta>0.
  \end{array}
 \right.  
\end{equation*}
If $w(s;\beta)$ has a zero in $(0,\infty)$, $\zeta(\beta)$ denotes the first 
zero of $w(s;\beta)$ in $(0,\infty)$, that is,
\begin{equation*}
 w(s;\beta)>0 \ \textup{for} \ 0<s<\zeta(\beta), \quad w(\zeta(\beta);\beta)=0.
\end{equation*}
Since $w'(\zeta(\beta);\beta)<0$, the implicit function theorem shows that 
$\zeta\in C^1$ near $\beta$, provided $\zeta(\beta)$ exists.

\begin{lemma}\label{zeta(A)}
 Let $p>1$, $h\in C[0,\infty)$ and $h(s)>0$ for $s>0$.
 For each $\beta_0>0$, there exists $A>\beta_0$ such that 
 $w(s;A)$ has a zero in $(0,\infty)$, that is, $\zeta(A)$ exists.
\end{lemma}

\begin{proof}
 Let $\beta_0>0$ and $M>0$ be a constant 
 as in the proof of Lemma \ref{w'(0)-->0}.
 Let $c>0$ satisfy
 \begin{equation*}
  \left(\frac{p+1}{M}\right)^\frac{1}{p-1}c^{-\frac{p+1}{p-1}} > \beta_0.
 \end{equation*}
 Let $w_c$ be a least energy solution of \eqref{bvpw}. 
 We set $A:=w_c'(0)$.
 Then $w(s;A)\equiv w_c(s)$ and $\zeta(A)=c$.
 We show $A>\beta_0$. 
 Let $m \in (0,c)$ satisfy $w_c'(m)=0$ and $\max_{s\in[0,c]}w_c(s)=w_c(m)$.
 Since $w_c$ is concave on $[0,m]$, we have
 \begin{equation*}
  w_c(s) \le w_c'(0)s=As, \quad s \in [0,c].
 \end{equation*}
 We observe that
 \begin{align*}
  A=w_c'(0) = - \int_0^m w_c''(s) ds
          = \int_0^m h(s) [w_c(s)]^p ds 
        & \le \int_0^c M A^p s^p ds \\
        & = \frac{MA^p}{p+1} c^{p+1}.
 \end{align*}
 Therefore,
 \begin{equation*}
  A \ge \left(\frac{p+1}{M}\right)^\frac{1}{p-1}c^{-\frac{p+1}{p-1}} > \beta_0.
 \end{equation*}
 This completes the proof.
\end{proof}

\begin{lemma}\label{existlargesol}
 Let $p>1$, $h\in C[0,\infty)$ and $h(s)>0$ for $s>0$.
 Assume that equation \eqref{w''+hw^p=0} has a solution $w_0$ 
 for which $w_0(0)=0$, $w_0(s)>0$ for $s>0$.
 Then, for all sufficiently large $c>0$, problem \eqref{bvpw} has a solution 
 $w_1$ satisfying $w_1'(0)>w_0'(0)$. 
\end{lemma}

\begin{proof}
 Let $\beta_0=w_0'(0)$. 
 Let $A$ be the constant as in Lemma \ref{zeta(A)}.
 By a general theory on the continuous dependence of solutions on initial 
 values, there exists $\overline{\beta_0} \in [\beta_0,A)$ such that 
 $\zeta(\beta)$ exists for each $\beta\in(\overline{\beta_0},A]$ and 
 $\zeta(\beta)\to\infty$ as $\beta\downarrow\overline{\beta_0}$.
\end{proof}

\begin{lemma}\label{existsmallandmidsol}
 Let $p>1$, $h\in C[0,\infty)$ and $h(s)>0$ for $s>0$.
 Assume that $h$ is bounded on $[0,\infty)$ and satisfy \eqref{s^sigmah(s)>0}
 and that equation \eqref{w''+hw^p=0} has a solution $w_0$ 
 for which $w_0(0)=0$, $w_0(s)>0$ for $s>0$.
 Then, for all sufficiently large $c>0$, problem \eqref{bvpw} has two solutions 
 $w_2$ and $w_3$ satisfying $0<w_3'(0)<w_2'(0)<w_0'(0)$. 
\end{lemma}

\begin{proof} 
 Let $\beta_0=w_0'(0)$. 
 By Lemma \ref{w'(0)-->0}, there exists $B\in(0,\beta_0)$ such that $\zeta(B)$ 
 exists.
 In the same way as in the proof of Lemma \ref{existlargesol}, 
 there exists $\underline{\beta_0} \in (B,\beta_0]$ such that 
 $\zeta(\beta)$ exists for each $\beta\in[B,\underline{\beta_0})$ and 
 $\zeta(\beta)\to\infty$ as $\beta\uparrow\underline{\beta_0}$.
 Let $c>0$ be arbitrary large.
 Then $\zeta(\beta_2)=c$ for some $\beta_2\in[B,\underline{\beta_0})$. 
 In addition, by Lemma \ref{w'(0)-->0}, there exists $\beta_3\in(0,B)$ 
 such that $\zeta(\beta_3)=c$.
 Consequently, $w(s;\beta_2)$ and $w(s;\beta_3)$ are solutions of \eqref{bvpw} 
 satisfying $0<w'(0;\beta_3)<B\le w'(0;\beta_2)<\underline{\beta_0}\le\beta_0$.  
\end{proof}

Theorem \ref{exist3solsODE} follows from Lemmas \ref{existlargesol} and 
\ref{existsmallandmidsol} immediately.

By the same arguments as in the proofs of Lemmas \ref{existlargesol} and 
\ref{existsmallandmidsol}, we can obtain the following result.

\begin{lemma}\label{exist2sols}
 Let $p>1$, $h\in C[0,\infty)$ and $h(s)>0$ for $s>0$.
 Suppose that equation \eqref{w''+hw^p=0} has two solutions $\underline{w}$ and 
 $\overline{w}$ for which $\underline{w}(0)=\overline{w}(0)=0$,
 $0<\underline{w}'(0)<\overline{w}'(0)$, 
 $\underline{w}(s)>0$ and $\overline{w}(s)>0$ for $s>0$.
 Assume moreover that equation \eqref{w''+hw^p=0} has a solution $w$ such that 
 $w(0)=0$, $0<\underline{w}'(0)<w'(0)<\overline{w}'(0)$, $w(s)>0$ for $0<s<s_0$,
 and $w(s_0)=0$ for some $s_0>0$.
 Then, for all sufficiently large $c>0$, problem \eqref{bvpw} has two solutions 
 $w_1$ and $w_2$ satisfying 
 $\underline{w}'(0)<w_1'(0)<w'(0)<w_2'(0)<\overline{w}'(0)$. 
\end{lemma}

\section{Proof of Theorem \ref{existksols}}

We will prove Theorem \ref{existksols}.
Throughout this section, we always assume that 
$N\in\mathbb{N}$, $N\ge3$, $(N+2)/(N-2)<p<p_{JL}$, $\lambda>0$, 
and let $r_*$ be a constants as in Theorem C.
Let $k\ge 2$ be an integer.
From Theorem C, there exists $\delta_1\in(0,r_*)$ such that 
problem \eqref{B} has at least $k$ positive solutions if $|b-r_*|<\delta_1$.

Let $u(r;\alpha)$ denote a unique solution of the initial value problem
\begin{equation}\label{IVP1}
 \left\{
  \begin{array}{l}
   u'' + \dfrac{N-1}{r} u' + \lambda u + |u|^{p-1}u = 0, \quad r >0, 
    \\[2ex]
   u(0)=\alpha>0, \quad u'(0)=0.
  \end{array}
 \right.
\end{equation}

\begin{proposition}\label{existz(alpha)}
 For each $\alpha>0$, $u(r;\alpha)$ has a zero in $(0,\infty)$. 
\end{proposition}

\begin{proof}
 The solution $\phi$ of \eqref{phi} has a zero $j_1\in (0,\infty)$, since 
  \begin{equation*}
   \phi(r)=cr^{-\frac{N-2}{2}}J_{\frac{N-2}{2}}(\sqrt{\lambda}r),
  \end{equation*}
 where $J_\mu$ is the Bessel function of the first kind of order $\mu$ and 
 $c$ is a specific positive constant.
 Assume to the contrary that $u:=u(r;\alpha)>0$ for $r\ge 0$.
 By the same argument as in the proof of Lemma \ref{phi>0}, 
 we lead a contradiction.
 Namely,
 \begin{align*}
  j_1^{N-1} \phi'(j_1) u(j_1) 
   & = \int_0^{j_1} 
   \left[ r^{N-1} (\phi'(r)u(r)-\phi(r)u'(r) ) \right]' dr \\
   & = \int_0^{j_1} r^{N-1} |u(r)|^{p-1} u(r) \phi(r) dr > 0,
 \end{align*}
 which contradicts the fact that $\phi'(j_1)<0$ and $u(j_1)>0$.
\end{proof}

Let $z(\alpha)$ denote the the smallest zero of $u(r;\alpha)$ in $(0,\infty)$.
Since $u(z(\alpha);\alpha)=0$ and $u'(z(\alpha);\alpha)<0$, 
the implicit function theorem shows that $z\in C^1(0,\infty)$.

Merle-Peletier \cite[Theorem 1.1]{MP} proved that the equation
\begin{equation*}
 u'' + \dfrac{N-1}{r} u' + \lambda u + |u|^{p-1}u = 0
\end{equation*}
has a unique radial singular solution $u_*$, and it satisfies
\begin{equation}\label{u*}
 u_*(r) = A r^{-\frac{2}{p-1}}(1+o(1)) \quad \textup{as} \ r \to 0,
\end{equation}
where
\begin{equation*}
 A = \left[ \frac{2}{p-1}\left( N-2-\frac{2}{p-1} \right)\right]^{1/(p-1)}.
\end{equation*}
See also Miyamoto and Naito \cite[Theorem 1.1]{MN2020}.

\begin{proposition}\label{r*}
 There exists $r_*>0$ such that $u_*(r_*)=0$ and $u_*(r)>0$ for $0<r<r_*$.
\end{proposition}

\begin{proof}
 This result must be known.  
 We give a proof for the reader's convenience.
 We note that $u_*$ satisfies 
 \begin{equation*}
  0\le -u_*'(r) \le C r^{-\frac{2}{p-1}-1}, \quad 0<r\le r_0 
 \end{equation*}
 for some $r_0>0$.
 See, for example, Miyamoto and Naito \cite[Lemma 2.1]{MN2020}.
 Hence, since $p>(N+2)/(N-2)$, $u_*$ satisfies 
 \begin{equation*}
  \lim_{r\to0^+} r^{N-1} u_*'(r) = 0.
 \end{equation*}
 From \eqref{u*}, it follows that $\lim_{r\to0^+} r^{N-2} u_*(r) = 0$,
 which implies 
 \begin{equation*}
  \lim_{r\to0^+} r^{N-1} u_*(r) = 0.
 \end{equation*}
 By the same argument as in the proof of Proposition \ref{existz(alpha)},
 we conclude that $u_*(r)$ has a zero in $(0,\infty)$.
\end{proof}

Hereafter, let $r_{*}$ be the constant as in Proposition \ref{r*}.

\begin{lemma}\label{alpha_i}
 There exists a sequence $\{\alpha_i\}$ such that
 \begin{gather*}
  0<\alpha_1<\alpha_2<\cdots<\alpha_i<\alpha_{i+1}<\cdots,
  \quad \alpha_i \to \infty, \\
  0<z(\alpha_{2j-1})<r_*, \ \ r_*<z(\alpha_{2j})<\infty, \quad j \in \mathbb{N}.
 \end{gather*}
\end{lemma}

\begin{proof}
First we show that $z(\alpha)\to r_*$ as $\alpha\to\infty$.
By Miyamoto and Naito \cite[Theorem 1.1]{MN2020}, we have
\begin{align*}
 u(r;\alpha) \to u_*(r) \quad \textup{in} \ C_{loc}^2(0,r_0] \ 
 \textup{as} \ \alpha \to \infty
\end{align*}
for some $r_0\in(0,r_*)$.
Thus, $(u(r_0;\alpha),u'(r_0;\alpha))$ is very close to \linebreak
$(u_*(r_0),u_*'(r_0))$ when $\alpha$ is sufficiently large.
By a general theory on the continuous dependence of solutions on initial values,
we conclude that 
\begin{align*}
 u(r;\alpha) \to u_*(r) \quad \textup{in} \ C^1_{loc}[r_0,r_*) \ 
 \textup{as} \ \alpha \to \infty.
\end{align*}
Let $\varepsilon>0$. 
There exists $\alpha_\varepsilon>0$ such that
$u(r_*-\varepsilon;\alpha)\ge u_*(r_*-\varepsilon)/2>0$ 
for $\alpha\ge\alpha_\varepsilon$. 
Since $u(r;\alpha)$ is decreasing with respect to $r\in(0,z(\alpha))$,
we have $u(r;\alpha)\ge u(r_*-\varepsilon;\alpha)>0$ for 
$r\in[0,r_*-\varepsilon]$ when $\alpha \ge \alpha_\varepsilon$,
which means that $z(\alpha)\ge r_*-\varepsilon$ for 
$\alpha \ge \alpha_\varepsilon$.
Therefore, $\liminf_{\alpha\to\infty}z(\alpha)\ge r_*$.

It is sufficient to prove that $\limsup_{\alpha\to\infty}z(\alpha)\le r_*$.
Assume that $\limsup_{\alpha\to\infty}z(\alpha)>r_*$.
Then there exists $\{a_i\}_{i=1}^\infty$ such that $z(a_i)>r_*$ 
for $i\ge 1$ and $\lim_{i\to\infty}z(a_i)\in (r_*,\infty]$.
Since $(r^{N-1}u'(r;\alpha))' \le 0$ for $r\in[0,z(\alpha)]$,
for each $r_1\in(0,z(\alpha))$, we get
\begin{equation*}
 r^{N-1}u'(r;a_i) \le r_*^{N-1} u'(r_*;a_i), 
 \quad r \in [r_*,z(a_i)].
\end{equation*}
Multiplying the both sides by $r^{1-N}$ and integrating it over $[r_*,z(a_i)]$, 
we have
\begin{equation*}
 0 = u(z(a_i);a_i) \le u(r_*;a_i)
 - \frac{r_*^{N-1} u'(r_*;a_i)}{N-2} (z(a_i)^{2-N}-r_*^{2-N}).
\end{equation*}
%Since $u'(r_*;a_i)<0$ for sufficiently large $i$, we have
%\begin{equation*}
% z(a_i)^{2-N} \ge r_*^{2-N} +\frac{(N-2)u(r_*;a_i)}{r_*^{N-1} u'(r_*;a_i)}.
%\end{equation*}
This contradicts the fact that 
$\lim_{i\to\infty}z(a_i)\in(r_*,\infty]$,
since $u(r_*;a_i)\to u_*(r_*)=0$ and $u(r_*;a_i)'\to u_*'(r_*)\in(-\infty,0)$.
Consequently, $z(\alpha)\to r_*$ as $\alpha\to\infty$.

Next we claim that $\#\{r\in(0,\min\{r_*,z(\alpha)\}): u(r;\alpha)=u_*(r)\}$ 
is finite for each fixed %$r_1>0$ and
$\alpha>0$.
We suppose that there exists 
$\{c_j\}_{j=1}^\infty\subset(0,\min\{r_*,z(\alpha)\})$ such that 
$u(c_j;\alpha)=u_*(c_j)$ for $j\in\mathbb{N}$.
Then $\{c_j\}$ contains a subsequence $\{c_{j_k}\}$ such that $c_{j_k}\to c$
for some $c \in [0,\min\{r_*,z(\alpha)\}]$.
Since $u(0;\alpha)=\alpha$ and $u_*(r)\to\infty$ as $r\to0$, we find $c\ne0$ and
$u(c;\alpha)=u_*(c)$.
Rolle's theorem implies $u'(c;\alpha)=u_*'(c)$.
By the uniqueness of solutions to the initial value problem, we get
$u(r;\alpha)\equiv u_*(r)$ for $r\in(0,\min\{r_*,z(\alpha)\})$, 
which is a contradiction.  

Assume, to the contrary, that there exists $A>0$ such that
either $0<z(\alpha)\le r_*$ for every $\alpha\ge A$ or 
$z(\alpha)\ge r_*$ for every $\alpha\ge A$.
We note that $\#\{r\in(0,\min\{r_*,z(\alpha)\}): u(r;A)=u_*(r)\}$ 
is a nonnegative integer.
From Miyamoto \cite[Theorem A]{M2018}, it follows that
$\#\{r\in(0,r_*/2) : u(r;\alpha)=u_*(r)\}\to\infty$ as $\alpha\to\infty$.
Recalling $z(\alpha)\to r_*$ as $\alpha\to\infty$, we see that
$\#\{r\in(0,\min\{r_*,z(\alpha)\}) : u(r;\alpha)=u_*(r)\}\to\infty$ 
as $\alpha\to\infty$.
%(In \cite{M2018}, only nonnegative solutions have been considered.)
By the continuity of $u(r;\alpha)$ with respect to $(r,\alpha)$ 
and the assumption that either $0<z(\alpha)\le r_*$ for every $\alpha\ge A$ or 
$z(\alpha)\ge r_*$ for every $\alpha\ge A$,
there exist $A_1\ge A$ and $r_1 \in (0,\min\{r_*,z(A_1)\}]$ such that 
$u(r_1;A_1)=u_*(r_1)$ and $u'(r_1;A_1)=u_*'(r_1)$.
Then $u(r;A_1)\equiv u_*(r)$ for $r\in(0,\min\{r_*,z(A_1)\})$,
by the uniqueness of solutions to the initial value problem.
This is a contradiction.
\end{proof}

Let $k$ be a positive integer.
We set 
\begin{equation*}
 \rho_k:=\frac{1}{2}\min\{|z(\alpha_i)-r_*| : 1 \le i \le 2k+1 \}>0.
\end{equation*}
Let $b\in[r_*-\rho_k,r_*+\rho_k]$  be arbitrary fixed.
Then we have 
\begin{align*}
 z(\alpha_{2j-1}) & \le r_* - 2\rho_k < b, \quad 1 \le j \le k+1, \\
   z(\alpha_{2j}) & \ge r_* + 2\rho_k > b, \quad 1 \le j \le k.
\end{align*}
Hence, there exists $\{\beta_i\}_{i=1}^{2k}$ such that
\begin{gather*}
 z(\beta_i)=b, \ \ \quad i=1,2,\cdots,2k, \\
 \alpha_i<\beta_i<\alpha_{i+1}, \quad i=1,2,\cdots,2k.
\end{gather*}
In addition, if $1\le l < m \le 2k$, then there 
exists $\eta_0 \in (\beta_l,\beta_m)$ such that $z(\eta_0)\ne b$.

Let $v(r;\gamma)$ denote a unique solution of initial value problem 
\begin{equation}\label{IVP2}
 \left\{
  \begin{array}{l}
   v'' + \dfrac{N-1}{r} v' + \lambda v + |v|^{p-1}v = 0, \quad 0 <r<b, 
    \\[2ex]
   v(b)=0, \ v'(b)=-\gamma<0.
  \end{array}
 \right.
\end{equation}
Since $u(b;\beta_i)=0$, if $-\gamma=u'(b;\beta_i)$, then
$v(r;\gamma)\equiv u(r;\beta_i)$.
Thus, there exists $\{\gamma_m\}_{m=1}^{2k}$ such that
\begin{gather*}
 0<\gamma_1<\gamma_2<\cdots<\gamma_{2k}, \\
 v(r;\gamma_m) \equiv u(r;\beta_{i(m)}) \quad \mbox{for\ some\ } 
 i(m) \in \{1,2,\cdots,2k\}.
\end{gather*} 
Moreover, if $m\ne m'$, then $i(m)\ne i(m')$.

\begin{lemma}\label{SorR}
 Let $\gamma>0$.
 Assume $v(r;\gamma)>0$ for $r\in(0,b)$.
 Then $v'(r;\gamma)<0$ for $r\in(0,b)$ and one of the following \textup{(i)} 
 and \textup{(ii)} holds\textup{:}
 \begin{enumerate}
  \item $\lim_{r\to0}v(r;\gamma)=\infty$\textup{;}
  \item $0<v(0;\gamma)<\infty$ and $v'(0;\gamma)=0$.
 \end{enumerate}
\end{lemma}

\begin{proof}
 We set $v(r):=v(r;\gamma)$.
 First we prove that $v'(r)<0$ for $r\in(0,b)$.
 We see that $V(t):=v(t^{-1/(N-2)})$ satisfies $V(t)>0$ for $t>b^{-(N-2)}$ and 
 is a solution of \eqref{EqU}.
 Then $V''(t)<0$ for $t>b^{-(N-2)}$, which implies that 
 $V'(t)$ is decreasing in $t>b^{-(N-2)}$.
 Thus, there exists $L:=\lim_{t\to\infty}V'(t)\in[-\infty,\infty)$.
 Recalling that $V(t)>0$ for $t>b^{-(N-2)}$, we get $L\in[0,\infty)$.
 Therefore, $V'(t)>0$ for $t>b^{-(N-2)}$, which means
 $v'(r)<0$ for $r\in(0,b)$.

 Thus, either $\lim_{r\to0}v(r)=\infty$ or $0<v(0)<\infty$.
 Hereafter, we suppose $0<v(0)<\infty$.
 Then $L=0$.
 We observe that
 \begin{align*}
  \lim_{r\to0} r^{N-1} v'(r) = \lim_{t\to\infty}(-(N-2)V'(t)) = -(N-2)L=0. 
 \end{align*}
 Integrating $-(r^{N-1}v'(r))'=r^{N-1}(\lambda v(r)+(v(r))^p)$ over 
 $[\varepsilon,r]$ and letting $\varepsilon\to0^+$, we get
 \begin{equation}\label{vdot}
  -r^{N-1}v'(r) = \int_0^r s^{N-1} (\lambda v(s)+(v(s))^p) ds, \quad 0<r<b.
 \end{equation}
 Therefore,
 \begin{align*}
  |v'(r)| & \le \frac{1}{r^{N-1}} \int_0^r s^{N-1} (\lambda v(s)+(v(s))^p) ds \\
  & \le \frac{1}{r^{N-1}} (\lambda v(0)+(v(0))^p) \int_0^r s^{N-1} ds \\
  & = \frac{1}{N} (\lambda v(0)+(v(0))^p) r, \quad 0<r<b. 
 \end{align*}
 Letting $r\to 0$, we have $v'(0)=0$.
\end{proof}

\begin{lemma}\label{existszero}
 Let $m\in\{1,2,\cdots,2k\}$.
 Suppose that $-u_*'(b)\not\in(\gamma_m,\gamma_{m+1})$ if $b=r_*$. 
 Then there exists $\zeta_m\in(\gamma_m,\gamma_{m+1})$ such that 
 $v(r;\zeta_m)$ has a zero in $(0,b)$.
\end{lemma}

\begin{proof}
 We assume that there is no such a $\zeta_m$.
 Then, for each $\gamma\in[\gamma_m,\gamma_{m+1}]$, 
 we have $v(r;\gamma)>0$ for $0<r<b$.
 By the uniqueness of singular solution $u_*$, 
 Lemma \ref{SorR} implies that $v(0;\gamma)\in(0,\infty)$ and 
 $v'(0;\gamma)=0$ for every $\gamma\in[\gamma_m,\gamma_{m+1}]$.

 Now we prove that $v(0;\gamma)$ is continuous with respect to 
 $\gamma\in[\gamma_m,\gamma_{m+1}]$.
 Let $\xi_0\in[\gamma_m,\gamma_{m+1}]$ and let 
 $\{\xi_n\}\subset[\gamma_m,\gamma_{m+1}]$ be an arbitrary sequence converging 
 to $\xi_0$.
 We set $\tau_n:=v(0;\xi_n)$. 
 Then $v(r;\xi_n)\equiv u(r;\tau_n)$.
 
 We claim that $\{\tau_n\}$ is bounded.
 Assume that $\{\tau_n\}$ is unbounded.
 Then a subsequence $\{\tau_{n_j}\}$ satisfies $\tau_{n_j}\to\infty$.
 By Miyamoto and Naito \cite[Theorem 1.1]{MN2020}, we find that
 \begin{align*}
  u(r;\tau_{n_j}) \to u_*(r) \quad \textup{in} \ C_{loc}^1(0,b].
 \end{align*}
 (Recall that the argument as in the proof of Lemma \label{alpha_i}.)
 This means that $b=r_*$.
 Indeed, since $u(r;\tau_{n_j})>0$ and $u'(r;\tau_{n_j})<0$ from \eqref{vdot} 
 for $0<r<b$ and 
 $u(b;\tau_{n_j})=0$, we get $u_*(r)\ge 0$ and $u_*'(r)\le 0$ for $0<r<b$ and 
 $u_*(b)=0$.
 If $u_*(r)$ has a zero $r_1\in(0,b)$, then $u_*(r)=0$ on $[r_1,b]$,
 which implies $u_*'(b)=0$.
 By the uniqueness of solutions to initial value problems, we conclude that
 $u_*(r)\equiv 0$, which is a contradiction.
 Hence, $u_*(r)>0$ for $0<r<b$ and we have $b=r_*$.
 Since $v(r;\xi_n)\equiv u(r;\tau_n)$ and 
 $\xi_n\to\xi_0\in[\gamma_m,\gamma_{m+1}]$, we get $v(r;\xi_0)\equiv u_*(r)$.
 However, this contradicts the fact that $v(0;\gamma)\in(0,\infty)$ 
 for every $\gamma\in[\gamma_m,\gamma_{m+1}]$.
 Consequently, $\{\tau_n\}$ is bounded as claimed.

 Integrating $-(r^{N-1}v'(r))'=r^{N-1}(\lambda v(r)+(v(r))^p)$ over $[0,s]$, 
 dividing both sides by $s^{N-1}$, and integrating it over $[0,r]$, we get 
 form $\lim_{r\rightarrow 0}v'(r;\xi_{n})=0$ by Lemma \ref{SorR}.
 \begin{align*}
  v(r;\xi_n) 
   & = v(0;\xi_n) - \int_0^r \frac{1}{s^{N-1}} \int_0^s t^{N-1}
   (\lambda v(t;\xi_n)+ v(t;\xi_n)^p) dt ds \\
   & = v(0;\xi_n) - \frac{1}{N-2} \int_0^r 
     t \left(1-\left(\frac{t}{r}\right)^{N-2}\right)
   (\lambda v(t;\xi_n)+ v(t;\xi_n)^p) dt
 \end{align*}
 for $r\in[0,b]$.
 This implies
 \begin{equation}
  0\le v(r;\xi_n) \le v(0;\xi_n)=\tau_n, \quad 0\le r \le b
   \label{U<tau_n}
 \end{equation}
 and
 \begin{equation*}
   v(0;\xi_n) = \frac{1}{N-2} \int_0^b 
     t \left(1-\left(\frac{t}{r}\right)^{N-2}\right)
   (\lambda v(t;\xi_n)+ v(t;\xi_n)^p) dt.
 \end{equation*}
 By the continuous dependence of solutions on initial values,
 for each $r_1\in(0,b)$, $v(r;\xi_n)$ converges to 
 $v(r;\xi_0)$ on $[r_1,b]$ uniformly.
 Hence, $v(r;\xi_n)$ converges to $v(r;\xi_0)$ for each fixed $r \in (0,b]$.
 Noting that $\{\tau_n\}$ is bounded and \eqref{U<tau_n}, 
 by Lebesgue's dominated convergence theorem, we have
 \begin{align*}
  \lim_{n\to\infty} v(0;\xi_n) 
   & = \lim_{n\to\infty} \frac{1}{N-2} \int_0^b 
     t \left(1-\left(\frac{t}{r}\right)^{N-2}\right)
   (\lambda v(t;\xi_n)+ v(t;\xi_n)^p) dt \\
   &  = \frac{1}{N-2} \int_0^b 
     t \left(1-\left(\frac{t}{r}\right)^{N-2}\right)
   (\lambda v(t;\xi_0)+ v(t;\xi_0)^p) dt \\
   &  = v(0;\xi_0).
 \end{align*}
 Consequently, $v(0;\gamma)$ is continuous with respect to 
 $\gamma\in[\gamma_m,\gamma_{m+1}]$.

 There exist $i(m)$ and $i(m+1)$ such that $i(m)\ne i(m+1)$,
 $v(r;\gamma_m)=u(r;\beta_{i(m)})$, $v(r;\gamma_{m+1})=u(r;\beta_{i(m+1)})$.
 Hereafter we suppose $\beta_{i(m)}<\beta_{i(m+1)}$. 
 The case $\beta_{i(m)}>\beta_{i(m+1)}$ can be treated similarly.
 Thus, we can take $\eta_0\in(\beta_{i(m)},\beta_{i(m+1)})$ for which
 $z(\eta_0)\ne b$.
 Since $v(0;\gamma)$ is continuous with respect to 
 $\gamma\in[\gamma_m,\gamma_{m+1}]$, $v(0;\gamma_m)=\beta_{i(m)}$ and
 $v(0;\gamma_{m+1})=\beta_{i(m+1)}$, there exists 
 $\gamma_0\in(\gamma_m,\gamma_{m+1})$ such that $v(0;\gamma_0)=\eta_0$.
 This means that $u(r;\eta_0)\equiv v(r;\gamma_0)$ and $z(\eta_0)=b$,
 which is a contradiction.
 Consequently, there exists $\zeta_m\in(\gamma_m,\gamma_{m+1})$ such that 
 $v(r;\zeta_m)$ has a zero in $(0,b)$.
\end{proof}

\begin{proof}[Proof of Theorem \ref{existksols}]
 We employ the transformation
 \begin{equation*}
  w(s;\gamma) = \frac{V(t;\gamma)}{\Phi(t)}, \quad 
   s = \int_{b^{-(N-2)}}^t \frac{1}{[\Phi(\tau)]^2} d\tau,
 \end{equation*}
 where $V(t;\gamma):=v(t^{-1/(N-2)};\gamma)$ and $\Phi$ is the function as in
 Section 2.
 Recalling Section 2, we conclude that $w=w(s;\gamma)$ satisfies
 \begin{equation*}
 \left\{
  \begin{array}{l}
   w'' + g(s)|w|^{p-1} w = 0, \quad 0<s<c, \\[1ex]
   w(0)=0, \quad w'(0)=\frac{b^{N-1}}{N-2}\gamma>0,
  \end{array}
 \right.  
 \end{equation*}
 where $g$ is the function defined by \eqref{g}.
 Let $m\in\{1,2,\cdots,2k\}$.
 Suppose that $-u_*'(b)\not\in(\gamma_m,\gamma_{m+1})$. 
 By Lemma \ref{existszero}, there exists $\zeta_m\in(\gamma_m,\gamma_{m+1})$ 
 such that $v(r;\zeta_m)$ has a zero in $(0,b)$.
 From Lemma \eqref{exist2sols}, there exists $a_m\in(0,b)$ such that
 if $0<a<a_m$, then problem \eqref{A} has two solutions $v_1$ and $v_2$ 
 such that $\gamma_m<-v_1'(b)<-v_2'(b)<\gamma_{m+1}$. 
 Since $m\in\{1,2,\cdots,2k\}$ and an $m'\in\{1,2,\cdots,2k\}$ may satisfy
 $-u_*'(b)\not\in(\gamma_{m'},\gamma_{m'+1})$, there are at least
 $2(2k-1)$ solutions of \eqref{A} when $0<a<a_0$ for some $a_0\in(0,b)$.
\end{proof}

\section{Proof of Theorem \ref{uniquenessresult}}
To prove Theorem \ref{uniquenessresult}, we outline the results from \cite{SW} .

\begin{proposition}[{\cite[Proposition 1]{SW}}]
Let $A\in C^{3}[a,b]\cap C[a,b]$.
For each solution $u$ of problem \eqref{A}, set 
\begin{align}\label{pohozaev}
&J(r;u):=\Bigl[
\frac{1}{2}A(r)u'(r)^2+\left(-\frac{1}{2}A_{r}(r)+\frac{N-1}{r}A(r)\right)u'(r)u(r)\\
&+\frac{1}{4}
\left(A''(r)-\frac{3(N-1))}{r}A'(r)+\frac{2N(N-1)}{r^{2}}A(r)\right)u(r)^{2}
+\frac{\lambda}{2}A(r)u(r)^2\nonumber\\
&+\frac{1}{p+1}A(r)u(r)^{p+1}\Bigr].\nonumber
\end{align}
Then, it holds
\begin{align}
\frac{dJ}{dr}(r;u)=G(r)u(r)^2+H(r)u(r)^{p+1},
\end{align}
where
\begin{align}
G(r)=&\frac{1}{4r^3}\Bigl[-4(N-1)(N+\lambda r^2)A(r)+r\left((N-1)(2N+3)\lambda r^2\right)A'(r)\nonumber\\
&-3r^2(N-1)A''(r)+r^3 A'''(r)\Bigr]\nonumber\\
H(r)=&\left[-\frac{N-1}{r}A(r)+\frac{p+3}{2(p+1)}A'(r)\right].\nonumber
\end{align}
\end{proposition}
The following is a key proposition for the proofs of Theorem \ref{uniquenessresult}.
\begin{proposition}[{\cite[Proposition 2]{SW}}]\label{keyproposition}
Let $A\in C^{3}(a,b)\cap C[a,b]$ satisfy 
\begin{enumerate}
 \item 
$A(a)=A(b)=0$,
\item
$G(r)\equiv 0$ on $[a,b]$,
\item
there exists $\kappa\in(a,b)$ such that 
$H(\kappa)=0$,
$H(r)>0$ on $(a,\kappa)$ and
$H(r)<0$ on $(\kappa,b)$,
\end{enumerate}
where $G$ and $H$ are the functions given in Proposition~\ref{pohozaev}.
Then problem \eqref{A} has a unique positive solution.
\end{proposition}
In the case of $N=3$, we can find that
\begin{align}\label{Ar}
A(r)=r^{2N-2}\cdot \frac{a\sin\sqrt{\lambda}(r-a)}{\sqrt{\lambda}r}\cdot\frac{b\sin\sqrt{\lambda}(b-r)}{\sqrt{\lambda}r}
\end{align}
satisfies conditions (i) and (ii) of Proposition \ref{keyproposition}; see proof of Lemmas 3 and 4 of \cite{SW}
(condition (i) is clearly satisfied). In this case, the form of the function $H$ is as follows:
\begin{lemma}\cite[Lemma 6]{SW}
\label{H-Z}
\[
H(r)= \frac{abr}{2\lambda(p+1)}Z(r),
\]
where
\begin{align*}
Z(r)
=&-\sqrt{\lambda}(p+3)r\sin\left(\sqrt{\lambda}(2r-a-b)\right)\\
&-2(p-1)
\sin\left(\sqrt{\lambda}(b-r)\right)
\sin\left(\sqrt{\lambda}(r-a)\right).
\end{align*}
\end{lemma}
Thus, all we have to do is to check condition (iii) of Proposition \ref{keyproposition}.
\begin{proof}[Proof of Theorem \ref{uniquenessresult}]
From Lemma~\ref{H-Z}, it is enough to show that 
there exists $\kappa\in (a,b)$ such that
$Z(r)>0$ on $[a,\kappa)$, 
$Z(\kappa)=0$
and
$Z(r)<0$ on 
$(\kappa,b]$.
From Proposition A.1 below it holds
\begin{align}\label{firsteigen}
 0<\sqrt{\lambda}(b-a)<\sqrt{\lambda_{1}(A(a,b))}(b-a)=\pi.
\end{align}
Thus, we obtain
\begin{equation}
\label{Gab} 
\left\{
\begin{aligned}
Z(a)=&a\sqrt{\lambda}(p+3)\sin\left(\sqrt{\lambda}(b-a)\right)>0,\\
Z(b)=&-b\sqrt{\lambda}(p+3)\sin\left(\sqrt{\lambda}(b-a)\right)<0.
\end{aligned}
\right.
\end{equation}
By simple calculations, 
we have
\[
 Z_{r}(r)=-2\lambda(p+3)r\cos\left(\sqrt{\lambda}(2r-a-b)\right)
+\sqrt{\lambda}(p-5)\sin\left(\sqrt{\lambda}(2r-a-b)\right).
\]
\begin{figure}[htb]
 \begin{center}
 \includegraphics[scale=0.85 ]{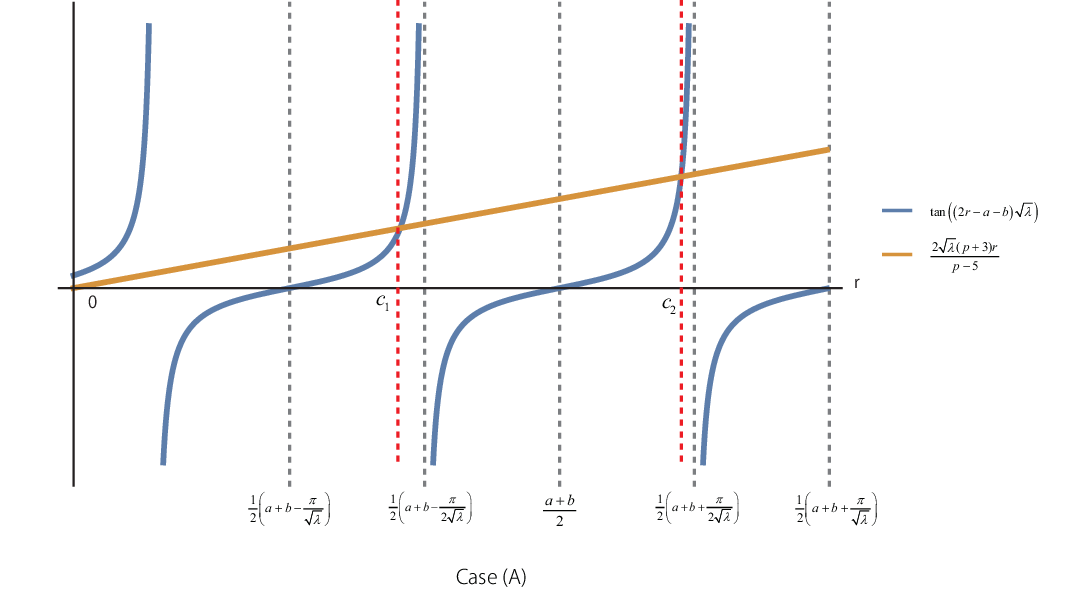}
\caption{ Graphs of $r\mapsto \tan\left(\sqrt{\lambda}(2r-a-b)\right)$ 
and $r\mapsto 2\sqrt{\lambda}(p+3)r/(p-5)$ for the Case (A).}
\label{fig2}
 \end{center}
\end{figure}
From \eqref{firsteigen},
we can find
\[
\frac{1}{2}\left(a+b-\frac{\pi}{\sqrt{\lambda}}\right)<a<b
<\frac{1}{2}\left(a+b+\frac{\pi}{\sqrt{\lambda}}\right)
\]
and
\begin{align}\label{rdom}
 -\pi < \sqrt{\lambda}(2r-a-b)
<\pi\Leftrightarrow\frac{1}{2}\left(a+b-\frac{\pi}{\sqrt{\lambda}}\right)
<r<\frac{1}{2}\left(a+b+\frac{\pi}{\sqrt{\lambda}}\right).
\end{align}

\vspace{0.5cm}
\noindent
Case (A):  $0\leq a+b-\frac{\pi}{\sqrt{\lambda}}$.\\
In this case, $Z$ must have exactly two critical points $r=c_1$(maximal point) and $r=c_2$(minimul point) when $r$ runs through the range \eqref{rdom}.
From this observation and \eqref{Gab} we can see that property (iii) of Proposition \ref{keyproposition} is satisfied; for more details, see proof of Theorem 1 of \cite{SW}.
Thus, we have proved the case $0\leq a+b-\frac{\pi}{\sqrt{\lambda}}$ and $b<a+\frac{\pi}{\sqrt{\lambda}}$.

\vspace{0.5cm}
\noindent
Case (B): $ a+b-\frac{\pi}{\sqrt{\lambda}}<0$ and $\frac{1}{2}(a+b-\frac{\pi}{2\sqrt{\lambda}})\leq a$.\\
We note $\frac{1}{2}(a+b-\frac{\pi}{2\sqrt{\lambda}})\leq a$ is equivalent to
\begin{align}\label{I1a}
 b\leq \frac{1}{2}(a+b+\frac{\pi}{2\sqrt{\lambda}})
\end{align}
Therefore, $Z$ is decreasing in $(a,(a+b)/2)$ and there exists at most one local minimum point of $Z$ in the interval 
$((a+b)/2,\frac{1}{2}(a+b-\frac{\pi}{2\sqrt{\lambda}}))$. Hence, noting \eqref{Gab}, property (iii) of Proposition \ref{keyproposition} is satisfied. 
The region of $(a,b)$ which satisfies 
\begin{align}\label{1111}
 a+b-\frac{\pi}{\sqrt{\lambda}}<0
\end{align} 
and \eqref{I1a} is labeled (B) in Figure \ref{fig:fig1}.
This proves the case.

Define an open intervals 
\begin{align*}
I_{1}&=\left( \frac{1}{2}(a+b-\frac{\pi}{\sqrt{\lambda}}),\frac{1}{2}(a+b-\frac{\pi}{2\sqrt{\lambda}})\right),\\
I_{2}&=\left(\frac{1}{2}(a+b), \frac{1}{2}(a+b+\frac{\pi}{2\sqrt{\lambda}})\right),\\
I_{3}&=\left(\frac{1}{2}(a+b+\frac{\pi}{2\sqrt{\lambda}}), \frac{1}{2}(a+b+\frac{\pi}{\sqrt{\lambda}})\right).
\end{align*}

\vspace{0.5cm}
\noindent
Case (C): $ a+b-\frac{\pi}{\sqrt{\lambda}}<0$ and $a\in I_{1}$.\\
Since 
\begin{align*}
\left(\tan\sqrt{\lambda}(2r-a-b)\right)'\vert_{r=\frac{1}{2}(a+b-\frac{\pi}{\sqrt{\lambda}})}=2\sqrt{\lambda}
\end{align*}
and $2\sqrt{\lambda}(p+3)/(p-5)>2\sqrt{\lambda}$, there exists a point $r_1 \in I_{1}$ such that derivatives of 
$\tan\sqrt{\lambda}(2r-a-b)$ and $2\sqrt{\lambda}(p+3)r/(p-5)$ are equal. Concretely, it holds
\begin{align*}
 \frac{2\sqrt{\lambda}}{\left(\cos\sqrt{\lambda}(2r_{1}-a-b)\right)^{2}}=\frac{2\sqrt{\lambda}(p+3)}{p-5}.
\end{align*}
Since $-\pi <\sqrt{\lambda}(2r_{1}-a-b)<-\pi/2$, we have 
\begin{align*}
 r_{1}=\frac{1}{2}\left(a+b-\frac{1}{\sqrt{\lambda}}\arccos\left(-\sqrt{\frac{p-5}{p+3}}\right)\right).
\end{align*}

\vspace{0.5cm}
\noindent
Case (C-i): $r_{1}\leq a$ and $a\in I_{1}$.
Note that  $r_{1}\leq a$ is equivalent to 
\begin{align}\label{r0a}
 b\leq a+\frac{1}{\sqrt{\lambda}}\arccos\left(-\sqrt{\frac{p-5}{p+3}}\right).
\end{align}
Also we note that $b\in I_{3}$. Hence, there are at most two critical points of $Z$, $c_{1}\in I_{1}$ and $c_{2}\in I_{2}$ ($c_{1}$ may not exists, but $c_{2}$ is necessarily exits). Since $c_{1}$ is maximal point and $c_{2}$ is minimal point, noting \eqref{Gab}, property (iii) of Proposition \ref{keyproposition} is satisfied.  
The region of $(a,b)$ which satisfies \eqref{1111} and \eqref{r0a} is labeled (C-i) in Figure \ref{fig:fig1}.
This proves the case.
\begin{figure}[htb]
 \begin{center}
 \includegraphics[scale=0.85 ]{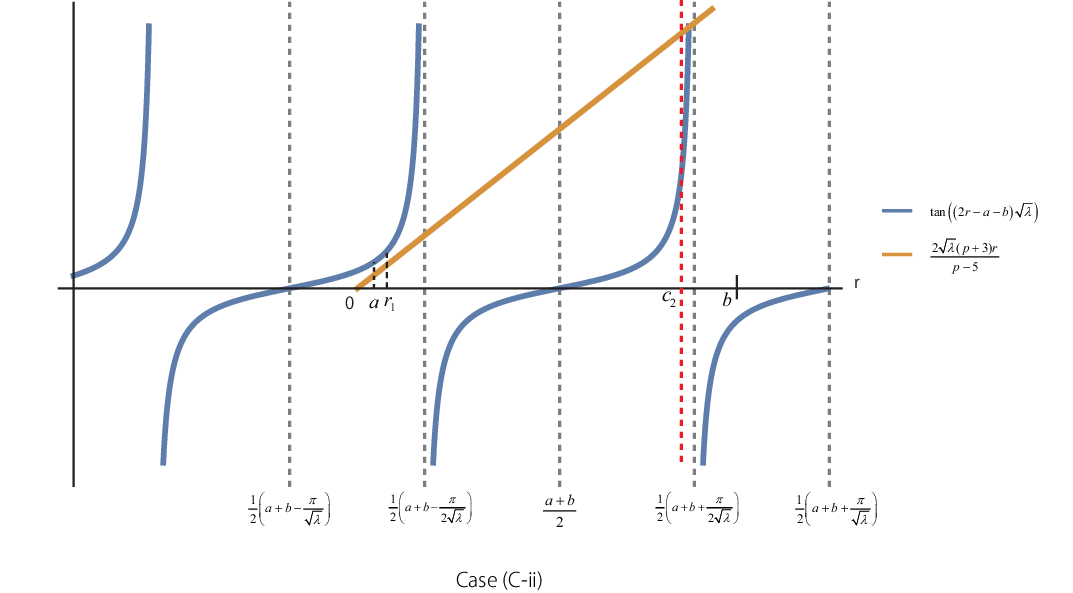}
\caption{ Graphs of $r\mapsto \tan\left(\sqrt{\lambda}(2r-a-b)\right)$ 
and $r\mapsto 2\sqrt{\lambda}(p+3)r/(p-5)$ for the Case (C-ii).}
\label{fig3}
 \end{center}
\end{figure}

\vspace{0.5cm}
\noindent
Case (C-ii): $a<r_{1}$ and $\tan\sqrt{\lambda}(2r_{1}-a-b)\geq \frac{2\sqrt{\lambda}(p+3)}{p-5}r_{1}$.\\
Note that $a<r_{1}$ is equivalent to 
\begin{align}\label{ar0}
 b> a+\frac{1}{\sqrt{\lambda}}\arccos\left(-\sqrt{\frac{p-5}{p+3}}\right).
\end{align} 
Also, $\tan\sqrt{\lambda}(2r_{1}-a-b)\geq \frac{2\sqrt{\lambda}(p+3)}{p-5}r_{1}$ is equivalent to
\begin{align}\label{zzz}
 b\leq -a +\frac{1}{\sqrt{\lambda}}\left(\frac{2\sqrt{2(p-5)}}{p+3}+\arccos\left(-\sqrt{\frac{p-5}{p+3}}\right)\right).
\end{align}
In this case, no critical point of $Z$ exists in $I_{1}$ and minimal point $c_{2}\in I_{2}$ exists; see Figure \ref{fig3}. Noting \eqref{Gab}, property (iii) of Proposition \ref{keyproposition} is satisfied.  
The region of $(a,b)$ which satisfies \eqref{1111}, \eqref{ar0} and \eqref{zzz} is labeled (C-ii) in Figure \ref{fig:fig1}.
This proves the case.
\begin{figure}[htb]
 \begin{center}
 \includegraphics[scale=0.85 ]{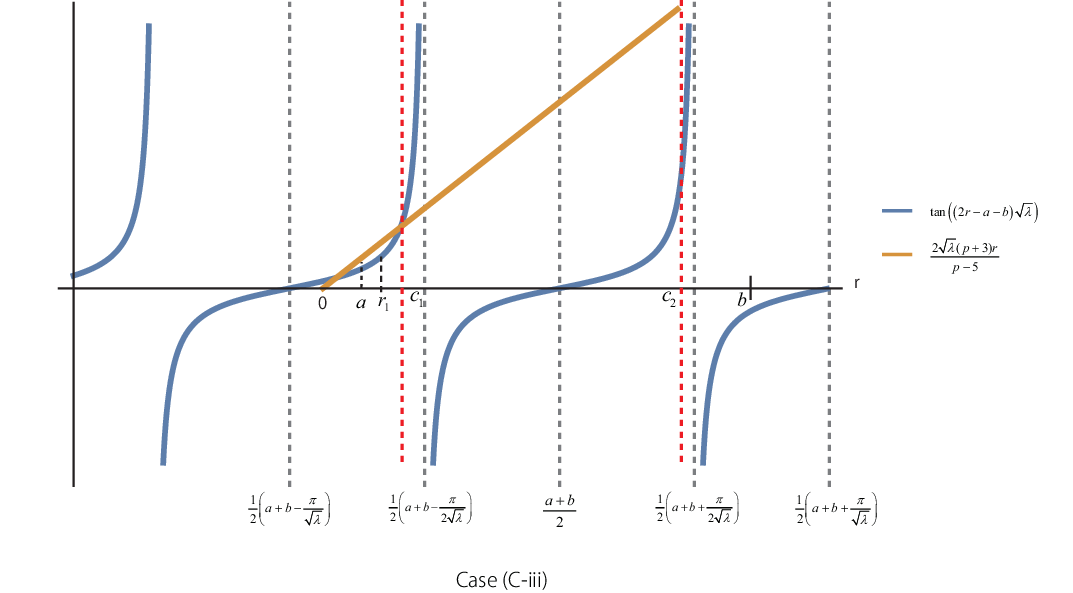}
\caption{ Graphs of $r\mapsto \tan\left(\sqrt{\lambda}(2r-a-b)\right)$ 
and $r\mapsto 2\sqrt{\lambda}(p+3)r/(p-5)$ for the Case (C-iii).}
\label{fig4}
 \end{center}
\end{figure}

\vspace{0.5cm}
\noindent
Case (C-iii): $a<r_{1}$, $\tan\sqrt{\lambda}(2r_{1}-a-b)< \frac{2\sqrt{\lambda}(p+3)}{p-5}r_{1}$ and 
$\tan\sqrt{\lambda}(a-b)\leq \frac{2\sqrt{\lambda}(p+3)}{p-5}a$.\\
As \eqref{zzz}, $\tan\sqrt{\lambda}(2r_{1}-a-b)< \frac{2\sqrt{\lambda}(p+3)}{p-5}r_{1}$ is equivalent to
\begin{align}\label{zzzz}
 b> -a +\frac{1}{\sqrt{\lambda}}\left(\frac{2\sqrt{2(p-5)}}{p+3}+\arccos\left(-\sqrt{\frac{p-5}{p+3}}\right)\right).
\end{align}
Noting that $-\pi<\sqrt{\lambda}(a-b)<-\pi/2$, we can see that $\tan\sqrt{\lambda}(a-b)\leq \frac{2\sqrt{\lambda}(p+3)}{p-5}a$ is equivalent to
\begin{align}\label{zzzzz}
b\geq a+\frac{1}{\sqrt{\lambda}}\left(\pi-\arctan\left(\frac{2\sqrt{\lambda}(p+3)}{p-5}a\right)\right).
\end{align}
In this case, maximal point $c_{1}$ of $Z$ exists in $I_{1}$ and minimal point $c_{2}\in I_{2}$ exists; see Figure \ref{fig4}. Noting \eqref{Gab}, property (iii) of Proposition \ref{keyproposition} is satisfied.  
The region of $(a,b)$ which satisfies \eqref{1111}, \eqref{ar0}, \eqref{zzz}  and \eqref{zzzzz} is labeled (C-iii) in Figure \ref{fig:fig1}.
This proves the case.
\end{proof}

\begin{proof}[Proof of Corollary \ref{<lam<}]
 Let $N=3$ and $p>5$.
 Then $\lambda_1(B(b))=\pi^2/b^{2}$.
 ($\phi_1(r)=r^{-1}\sin(\pi r/b)$ is an eigenfunction.)
 By Lemma A.2 below, if $\lambda\ge\lambda_1(B(b))=\pi^2/b^2$, then
 problem \eqref{B} has no positive solution. 
 Now we suppose 
 \begin{align}
 b<\dfrac{1}{\sqrt{\lambda}}\left(\frac{2\sqrt{2(p-5)}}{p+3}+\arccos\left(-\sqrt{\dfrac{p-5}{p+3}}\right)\right)=:\overline{b}
 \end{align}
 Then, Theorem \ref{uniquenessresult} implies that problem \eqref{A} has a unique 
 positive solution for each sufficiently small $a\in(0,b)$.
 %On the other hand, 
 Therefore, by Theorem \ref{exist3sols}, problem \eqref{B} has no positive solution.
 Indeed, if problem \eqref{B} has a positive solution, then 
 Theorem \ref{exist3sols} shows that problem \eqref{A} has at least three 
 positive solutions for all sufficiently small $a>0$,  
 which is a contradiction.

 Finally, we assume the case $b=\overline{b}$.
 In the same way as in the case $b<\overline{b}$, 
 it is sufficient to demonstrate that problem \eqref{A} has a unique positive 
 solution for each sufficiently small $a\in(0,b)$.
 By the same argument as in Case (C-ii) of the proof of Theorem 1.2,
 we conclude that if $a=0$ and $b=\overline{b}$, then there exists $c_2\in I_2$
 such that $Z_r(r)<0$ on $[0,r_1)$, $Z_r(r_1)=0$, $Z_r(r)<0$ on $(r_1,c_2)$,
 $Z_r(c_2)=0$, and $Z_r(r)>0$ on $(c_2,b]$.
 Since $Z(0)=0$ when $a=0$, we have $Z(r)<0$ for $r\in(0,b]$, especially 
 $Z(r_1)<0$.
 Let $a>0$ be sufficiently small.
 Then \eqref{zzzz} is satisfied.
 In a similar way as in Case (C-iii) of the proof of Theorem 1.2, 
 since $a>0$ is sufficiently small,
 there exist $c_0$, $c_1\in I_1$ and $c_2\in I_2$ such that
 $Z_r(r)<0$ on $[a,c_0)$, $Z_r(c_0)=0$, $Z_r(r)>0$ on $(c_0,c_1)$,
 $Z_r(c_1)=0$, $Z_r(r)<0$ on $(c_1,c_2)$, $Z_r(c_2)=0$, and 
 $Z_r(r)>0$ on $(c_2,b]$.
 Since $a>0$ is sufficiently small, $c_0$, $r_0$, $c_1$ are very close.
 Recalling $Z(r_1)<0$ when $a=0$, we can take a small $a_0>0$ such that 
 if $a\in(0,a_0)$, then $Z(c_0)<0$, $Z(r_0)<0$ and $Z(c_1)<0$,
 and hence $Z(r)>0$ on $[a,\kappa)$, 
 $Z(\kappa)=0$ and $Z(r)<0$ on $(\kappa,b]$ for some $\kappa\in(a,c_0)$.
 This means that problem \eqref{A} has a unique positive solution 
 for each sufficiently small $a\in(0,a_0)$.
\end{proof}

\bigskip

\begin{center}
 {\sc Appendix}
\end{center}

\bigskip

\noindent{\bf Proposition A.1.}
\it
If $\lambda\ge\lambda_1(A(a,b))$, then problem \eqref{A} has no positive 
solution.
\rm

\bigskip

\begin{proof}
 Let $\lambda\ge\lambda_1(A(a,b))$.
 We suppose that problem \eqref{A} has a positive solution $u$.
 Let $\phi_1$ be an eigenfunction corresponding to $\lambda_1(A(a,b))$.
 Then $\phi_1$ is radial from its positivity.
 Since $u$ and $\phi_1$ are solutions of 
 \begin{equation*}
  (r^{N-1} u')' + \lambda r^{N-1} (1 + \lambda^{-1}|u(r)|^{p-1}) u = 0;
  \quad u(a)=u(b)=0
 \end{equation*}
 and 
 \begin{equation}\label{phi1}
  (r^{N-1} \phi_1')' + \lambda_1(A(a,b)) r^{N-1} \phi_1 = 0; \quad
 \phi_1(a)=\phi_1(b)=0,
 \end{equation}
 respectively, the Strum comparison theorem implies that
 $u$ has a zero in $(a,b)$.
 This is a contradiction.
\end{proof}

Using the same argument as in the proofs of Proposition A.1 and 
Lemma \ref{phi>0}, we get the following result.

\bigskip

\noindent{\bf Proposition A.2.}
\it
If $\lambda\ge\lambda_1(B(b))$, then problem \eqref{B} has no positive solution.
\rm

\bigskip

\bigskip

\noindent{\bf Proposition A.3.}
\it
If $0<\lambda<\lambda_1(A(a,b))$, then problem \eqref{A} has a positive 
solution.
\rm

\bigskip

\begin{proof}
 We assume $0<\lambda<\lambda_1(A(a,b))$.
 Let $\phi_1$ be an eigenfunction corresponding to $\lambda_1(A(a,b))$ and
 let $\phi(r;\varepsilon)$ be a unique solution of the initial value problem
 \begin{equation}\label{phieps}
 \left\{
  \begin{array}{l}
   (r^{N-1}\phi')' +  \lambda r^{N-1} \phi = 0, \quad r>0, \\[2ex]
   \phi(a)=\varepsilon, \ \phi'(a)=1,
  \end{array}
 \right.
 \end{equation}
 where $\varepsilon\ge0$.
 Since $\phi_1$ is a nontrivial solution of \eqref{phi1},
 from the Strum comparison theorem, it follows that $\phi(r;0)>0$ for 
 $r\in(a,b]$. 
 By the continuous dependence of solutions on initial values,
 we conclude that $\phi(r;\varepsilon)>0$ for every sufficiently small 
 $\varepsilon>0$.
 We take such an $\varepsilon>0$ and set $\phi(r)=\phi(r;\varepsilon)$.
 By the transformations $\Phi(t)=\phi(t^{-1/(N-2)})$ and \eqref{Utow},
 problem \eqref{A} is transformed into \eqref{bvpw} for some $c>0$.
 Problem \eqref{bvpw} has a least energy solution for every $c>0$,
 and hence problem \eqref{A} has a positive solution.
\end{proof}

\begin{proof}[Proof of Proposition \ref{existpossol}]
 Proposition \ref{existpossol} follows immediately from Proposition A.1 and 
 A.3, 
\end{proof}

\bigskip

\noindent{\bf Lemma A.4.}
\it
Let $z(\alpha)$ be as in Section 3.
Then $z(\alpha)\to\sqrt{\lambda_1(B(1))/\lambda}$ as $\alpha\to0$.
\rm

\bigskip

\begin{proof}
 Let $u(r;\alpha)$ be a unique solution of the initial value problem 
 \eqref{IVP1}.
 We define the function $E(r)$ by
 \begin{equation*}
  E(r) = 
  \frac{1}{2} [u'(r;\alpha)]^2 + \frac{\lambda}{2} [u(r;\alpha)]^2
  + \frac{1}{p+1} |u(r;\alpha)|^{p+1}, \quad r\ge 0. 
 \end{equation*}
 Since
 \begin{equation*}
  E'(r) = - \frac{N-1}{r} [u'(r;\alpha)]^2 \le 0, \quad r>0,
 \end{equation*}
 we have
 \begin{equation*}
  [u(r;\alpha)]^2 \le \frac{2}{\lambda}E(r) \le \frac{2}{\lambda}E(0) 
   = \alpha^2 + \frac{2}{(p+1)\lambda}\alpha^{p+1}, \quad r \ge 0,
 \end{equation*}
 that is,
 \begin{equation}\tag{A.1}
  |u(r;\alpha)| \le 
   \left( 1 + \frac{2}{(p+1)\lambda}\alpha^{p-1} \right)^{1/2} \alpha, 
   \quad r \ge 0.
  \label{|u|<}
 \end{equation} 
 We set $\phi(r;\alpha)=\alpha^{-1}u(r;\alpha)$.
 Let $\phi_0$ be a solution of the initial value problem \eqref{phi}.
 We can prove that $\phi(r;\alpha)$ converges to a solution $\phi_0$ 
 as $\alpha\to0$ uniformly on $[0,R]$ for each fixed $R>0$.
 Indeed, since 
 \begin{equation*}
  \phi(r;\alpha) 
   = 1 -\frac{1}{N-2}\int_0^r s \left(1-\left(\frac{s}{r}\right)^{N-2}\right)
     (\lambda \phi(s;\alpha)+\alpha^{-1}|u(s;\alpha)|^{p-1}u(s;\alpha)) ds 
 \end{equation*}
 and
 \begin{equation*}
  \phi_0(r) 
   = 1 -\frac{1}{N-2}\int_0^r s \left(1-\left(\frac{s}{r}\right)^{N-2}\right)
     \lambda \phi_0(s) ds,
 \end{equation*}
 we see that $\psi(r):=\phi(r;\alpha)-\phi_0(r)$ satisfies
 \begin{align*}
  |\psi(r)| & \le \frac{\lambda}{N-2} \int_0^r s |\psi(s)| ds 
   + \frac{1}{\alpha(N-2)} \int_0^r s |u(s;\alpha)|^p ds.
 \end{align*}
 By using \eqref{|u|<}, Gronwall's inequality implies that $\psi(r)\to 0$
 as $\alpha\to0$ uniformly on $[0,R]$ for each fixed $R>0$.
 
 Therefore, $z(\alpha)$ converges to the smallest zero of $\phi_0(r)$ in 
 $(0,\infty)$ as $\alpha\to0$.
 Let $\phi_1$ be an eigenfunction corresponding to $\lambda_1(B(1))$
 satisfying $\phi_1(0)=1$.
 Then $\varphi_0(r):=\phi_1(\sqrt{\lambda/\lambda_1(B(1))}r)$ is a solution 
 of problem \eqref{phi}, which means $\phi_0(r)\equiv\varphi_0(r)$.
 We find that 
 \begin{equation*}
  \phi_0(\sqrt{\lambda_1(B(1))/\lambda})
 =\varphi_0(\sqrt{\lambda_1(B(1))/\lambda})=\phi_1(1)=0.  
 \end{equation*}
 Consequently, $z(\alpha)\to\sqrt{\lambda_1(B(1))/\lambda}$ as $\alpha\to0$.
\end{proof}

\begin{proof}[Proof of Theorem C]
 Miyamoto and Naito \cite[Theorem 1.1]{MN2020} proved
 \begin{align*}
  u(r;\alpha) \to u_*(r) \quad \textup{in} \ C_{loc}^1(0,b] \quad
  \textup{as} \ \alpha\to\infty,
 \end{align*}
 which implies that $z(\alpha)\to r_*$ as $\alpha\to\infty$, 
 where $r_*$ is the smallest zero of the singular solution $u_*$.
 By Lemma A.4, when $r_*<b<\sqrt{\lambda_1(B(1))/\lambda}$, there exists
 $\alpha\in(0,\infty)$ such that $z(\alpha)=b$, which proves (ii).

 Hereafter we assume that $(N+2)/(N-2)<p<p_{JL}$.
 By Dolbeault and Flores \cite[Theorem 1]{DF} or Guo and Wei 
 \cite[Theorem 1.1]{GW}, we see that, for each $k\in\mathbb{N}$, 
 problem \eqref{B} has at least $k$ positive solutions if
 $|b-r_*|$ is sufficiently small.
 Hence, if $|b-r_*|$ is sufficiently small, then there exists 
 $\{\rho_i\}_{i=1}^k$ such that $z(\rho_i)=b$.
 Moreover, 
 by Lemma A.4 and the fact that $z(\alpha)\to r_*$ as $\alpha\to\infty$,
 we conclude that
 $\underline{r}:=\min_{\alpha\in(0,\infty)}z(\alpha)$ exists
 and satisfies $0<\underline{r}<r_*$.
 Consequently, if $\underline{r}\le b<\sqrt{\lambda_1(B(1))/\lambda}$, then
 problem \eqref{B} has at least one positive solution.
\end{proof}

\begin{remark}
Let $N\in\mathbb{N}$, $N\ge3$ and $p>(N+2)/(N-2)$.
The existence of $\underline{r}:=\min_{\alpha\in(0,\infty)}z(\alpha)$
in the proof of Theorem C has been obtained in \cite{MPS} and \cite{PS}.
We can also show that $\underline{r}\ge r_0$, where $r_0$ is the constant 
as in Theorem A.
Indeed, if $\underline{r}<r_0$, then problem \eqref{B} with $b=\underline{r}$ 
has at least one positive solution and (iii) of Theorem A implies that 
problem \eqref{A} has at most one positive solution when $b=\underline{r}$ and 
$a\in(0,\underline{r})$.
This contradicts Theorem \ref{exist3sols}.
\end{remark}

{\bf Acknowledgements} \ 
The authors would like to thank the anonymous reviewer 
for the valuable comments.

\end{document}